\documentclass[12pt]{amsart}
\usepackage{bbm}
\usepackage{mathrsfs}
\usepackage[english]{babel}
\usepackage{wrapfig}

\usepackage[top=3cm,bottom=2cm,left=3cm,right=3cm,marginparwidth=1.75cm]{geometry}

\usepackage{amsmath,amssymb}
\usepackage{graphicx}
\usepackage{cite}
\usepackage{enumitem} 
\usepackage[all,cmtip]{xy}
\usepackage{tikz-cd}
\usetikzlibrary{intersections}
\usepackage{relsize}
\usepackage{faktor}
\usepackage{dsfont}
\usepackage{xcolor, import}
\usepackage{bm}
\usepackage{hyperref}
\usepackage{mathtools}
\hypersetup{
 citebordercolor={green!40!black},
 linkbordercolor={red!50!black},
 urlbordercolor={blue!70!black}
}

\makeatletter
\renewcommand{\section}{\@startsection%
{section}% name
{1}% level
{0mm}% indent
{1.5\bigskipamount}% beforeskip
{0.5\bigskipamount}% afterskip
{\centering\normalsize\sc}}% style

\renewcommand{\subsection}{\@startsection%
{subsection}% name
{2}% level
{0mm}% indent
{0.5\bigskipamount}% beforeskip
{0.5mm}% afterskip
{\normalsize\sc}}% style

\renewcommand{\paragraph}{\@startsection%
{paragraph}% name
{4}% level
{0mm}% indent
{\bigskipamount}% beforeskip
{0pt}% afterskip
{\normalsize\bf}}% style

\expandafter\let\expandafter\oldproof\csname\string\proof\endcsname
\let\oldendproof\endproof
\renewenvironment{proof}[1][\proofname]{%
  \oldproof[\slshape #1]%
}{\oldendproof}
\def\provedboxcontents#1{$\square$}
\newtheoremstyle{thm}{6pt plus 1pt minus 1pt}{6pt plus 1pt minus 1pt}{\slshape}{}{\scshape}{.}{5pt plus 1pt minus 1pt}{}
\newtheoremstyle{def}{6pt plus 1pt minus 1pt}{6pt plus 1pt minus 1pt}{}{}{\scshape}{.}{5pt plus 1pt minus 1pt}{}
\newtheoremstyle{rmk}{6pt plus 1pt minus 1pt}{6pt plus 1pt minus 1pt}{}{}{\scshape}{.}{5pt plus 1pt minus 1pt}{}
\newtheoremstyle{claim}{6pt plus 1pt minus 1pt}{6pt plus 1pt minus 1pt}{}{}{\slshape}{.}{5pt plus 1pt minus 1pt}{}

\theoremstyle{thm}
\newtheorem{newstatement}{newstatement}
\newtheorem{lemma}[newstatement]{Lemma}
\newtheorem{theorem}[newstatement]{Theorem}
\newtheorem*{theorem*}{Theorem 2}

\newtheorem*{stiefel-thm}{Stiefel's Parallelizability Theorem}

\theoremstyle{def}

\theoremstyle{rmk}
\newtheorem{remark}[newstatement]{Remark}

\newtheorem*{example*}{Example}

\theoremstyle{claim}

\renewcommand{\epsilon}{\varepsilon}
\renewcommand{\phi}{\varphi}

\newcommand{\R}{\mathbb{R}}

\newcommand{\Z}{\mathbb{Z}}

\newcommand{\defeq}{\vcentcolon=}

\let\emph\textsl
\title{On Dold-Whitney's parallelizability of 4-manifolds}

\author{Valentina Bais}
\address{Department of Mathematics, SISSA, Via Bonomea 265, 34136 Trieste, Italy.}
\email{vbais@sissa.it}

\begin{document}

%\normalpage

\begin{abstract}
We present a proof of the fact that a closed orientable 4-manifold is parallelizable if and only if its second Stiefel-Whitney class, first Pontryagin class and Euler characteristics vanish. This follows from a stronger result due to Dold and Whitney on the classification of oriented sphere bundles over a 4-complex. The contribution of this note is to outline in detail an argument which is essentially due to R. Kirby, using the classification of $SO(4)$-bundles over the 4-sphere by means of their Euler and first Pontryagin classes as a main tool.
\end{abstract}

\maketitle

Stiefel's Parallelizability Theorem asserts that any closed orientable 3-manifold has trivial tangent bundle. In addition to the original proof \cite{St1935}, one can find several elementary proofs of this fact, see for example Bais and Zuddas \cite{BZ}, Benedetti and Lisca \cite{BL2018}, Durst, Geiges, Gonzalo and Kegel \cite{DGGK2020}, Gonzalo \cite{Go1987}, Kirby \cite[Chapter VII]{Kirby}, Geiges \cite[Section 4.2]{Ge2008}, Fomenko and Matveev \cite[Section 9.4]{FM1997} and Whitehead \cite{Wh1961}. Motivated by this, one could ask whether it is possible to provide a simple proof of the following result. 

\begin{theorem}[Dold-Whitney]\label{main}
    A closed orientable 4-manifold $M$ is parallelizable if and only if its second Stiefel-Whitney class $w_2(M)$, first Pontryagin class $p_1(M)$ and Euler characteristic $\chi(M)$ are all vanishing.
\end{theorem}

Recall that an $n$-manifold $M$ is parallelizable if its tangent bundle $TM$ is trivial, i.e., if there is an isomorphism of vector bundles $TM \cong M \times \R^n$. One can easily verify that this is equivalent to ask for the existence of n vector fields on $M$ which are everywhere linearly independent. Such $n$-tuple is usually called a framing.

Theorem \ref{main} is due to Dold and Whitney's classification of $SO(4)$-bundles over closed, orientable 4-manifolds by means of their second Stiefel-Whitney class, first Pontryagin class and Euler characteristics, see \cite{Dold-Whitney}. A proof of this result can be also found in \cite{Thomas}. The argument we present is elementary and builds up on the discussion in Kirby's book \cite[Chapter VI]{Kirby}. It is worth to remark that it is still unknown which finitely presented fundamental groups can occur as the fundamental group of a closed parallelizable 4-manifold, see \cite{par}. Moreover, we recall that a closed orientable 4-manifold supports an Engel structure if and only if its tangent bundle is trivial, as shown in \cite{Vogel}.

We start by recalling the following argument on the classification of $SO(4)$-principal bundles over the 4-sphere. This can be found in \cite[Chapter VI]{Kirby}.

\begin{lemma}\label{lemma}
    $SO(4)$-principal bundles over $S^4$ are determined by their Euler and first Pontryagin classes. 
\end{lemma}
\begin{proof}
    Since $SO(4)$-principal bundles over $S^4$ are in one to one correspondence with homotopy classes of functions $f: S^3 \to SO(4)$, it is enough to show that the map
    \begin{equation}\label{isom}
        \Phi: \pi_3(SO(4)) \xrightarrow{(\chi, -(p_1+2 \chi)/4)} \Z \oplus \Z
    \end{equation}
    given by \[f \mapsto (\chi(f),-(p_1(f)+2 \chi(f))/4)\] is an isomorphism. Here $\chi(f)$ and $p_1(f)$ denote the Euler characteristics and the first Pontryagin class of the $SO(4)$-principal bundle associated to $f$ respectively. Note that we already know that $\pi_3(SO(4))\cong \pi_3(S^3 \times S^3) \cong \Z \oplus \Z$, since the universal cover of $SO(4)$ is $\text{Spin}(4) \cong S^3 \times S^3$.
    
    In the following, we view $S^3$ as the group of unit quaternions and work under the identifications \[\R^3\cong \langle i, j, k \rangle \subset \mathbb{H} \cong \R^4.\] 
    Let
    \begin{equation}
        \eta:S^3 \to SO(4)
    \end{equation}
    be the map defined by \[\eta(q)(q')\defeq q \cdot q'\] where $\cdot$ denotes the quaternionic multiplication. The homotopy class of this map determines the quaternionic Hopf bundle over $S^4$, which has Euler number $\chi(\eta)=1$. Its first Chern class $c_1(\eta) \in H^2(S^4; \Z)=0$ is necessarily trivial and this implies that the first Pontryagin class is \[p_1(\eta)=(c_1^2-2c_2)(\eta)=-2.\]

    We also define
    \begin{equation}
        \nu: S^3 \to SO(3) \subset SO(4)
    \end{equation}
    to be the map  \[\nu(q)(q')\defeq q \cdot q' \cdot q^{-1}\] and we compute the Euler and first Pontryagin class of the associated $SO(4)$-bundle. This can be done by noticing that the tangent bundle $\tau$ of $S^4$ is determined by the clutching function $(\eta)^2 \cdot (\nu)^{-1}$, where $\cdot$ and $(\ \cdot \ )^{-1}$ denote the point-wise multiplication and inversion in the Lie group $SO(4)$ respectively, see \cite[§23.6]{Steenrod}.  Since the characteristic class $p_1$ does not change under Whitney sums with trivial bundles, we also have that \[p_1(\tau)=p_1(\tau \oplus \epsilon^1)=p_1(\epsilon^5)=0.\] 
    Since $\chi(\tau)=2$, this implies that $\chi(\nu)=0$, $p_1(\nu)=-4$.

In particular, the map (\ref{isom}) is such that $\Phi(\eta)= (1,0)$ and $\Phi(\nu)=(0,1)$. Moreover, one can check that \[\Phi(f \cdot g)=\Phi(f)+\Phi(g)\] for every $f,g \in \pi_3(SO(4))$, where $\cdot$ denotes the point-wise multiplication in the Lie group $SO(4)$. This concludes the proof, since any surjective homomorphism $\Z \times \Z \rightarrow \Z \times \Z$ is a group isomorphism. 
\end{proof}

\begin{remark}
    The proof of Lemma \ref{lemma} also implies that, up to homotopy, the clutching functions of all possible $SO(4)$-bundles over the 4-sphere are obtained by multiplying copies of $\eta$, $\nu$ and their point-wise inverses.
\end{remark}

We are now ready to prove Theorem \ref{main}.
\begin{proof}
    Let $M$ be a closed orientable 4-manifold. It is straightforward to verify that the triviality of its tangent bundle implies the vanishing of $w_2(M)$, $p_1(M)$ and $\chi(M)$.

    Conversely, fix a handlebody decomposition of $M$ with just one 0-handle and one 4-handle. The vanishing of the second Stiefel-Whitney class implies the existence of a trivialization of the tangent bundle up to the 2-handles. Since $\pi_2(SO(4))=0$, one can always extend such trivialization over the 3-handles. The problem is now finding an extension over the remaining 4-handle $h^4$. In order to tackle this problem, we consider the change of basis map
    \[f: \partial(M \setminus h^4) \cong S^3 \to SO(4) \]
    between a fixed framing on $M \setminus h^4$ and the one of $\partial h^4$ that extends to the 4-ball. Note that, if $f$ is homotopic to the constant map, we can find a global framing for $M$, which is hence parallelizable.

    Consider the lift 
    \[ \widetilde f=(f_1,f_2): S^3 \to S^3 \times S^3 \cong \text{Spin} (4)\]
    of $f$ to the universal cover of $SO(4)$ and let $f_1,f_2$ be its two components. Since these are both continuous maps from the 3-sphere to itself, their degrees $\text{deg}(f_i)$ are defined for $i=1,2$. In particular, we will now show that 
    \begin{equation}\label{1}
        \text{deg} (f_1) - \text{deg} (f_2) = \chi(M)
    \end{equation}
    and 
    \begin{equation}\label{2}
        \text{deg} (f_1) + \text{deg} (f_2) = -p_1(M)/2,
    \end{equation}
    which will allow us to conclude our proof. Indeed, one has that $f$ is null-homotopic if and only if the degrees of both $f_1$ and $f_2$ are trivial, which by (\ref{1}) and (\ref{2}) happens only in the case in which both $\chi(M)$ and $p_1(M)$ vanish.
    
 We now show that equations \ref{1} and \ref{2} hold for the lifts \[(\eta_1, \eta_2), (\nu_1, \nu_2): S^3 \to S^3 \times S^3\] of $\eta$ and $\nu$ respectively. The conclusion will then follow from linearity. Indeed, it is an exercise to show that \[\text{deg}(a \cdot b)=\text{deg}(a)+\text{deg}(b)\] for any two continuous maps $a,b:S^3 \to S^3$ and, up to homotopy, every continuous function from the 3-sphere to $SO(4)$ can be expressed as a composition of copies of $\eta$, $\nu$ and the functions given by their point-wise quaternionic inverses.

 By definition, $\nu$ lifts to a map which is the identity of $S^3$ on both components and hence \[\text{deg}(\nu_1)-\text{deg}(\nu_2)=1-1=0=\chi(\nu)\] and
 \[\text{deg}(\nu_1) + \text{deg}(\nu_2)=1+1=2=-p_1(\nu)/2.\]

 On the other hand, $\eta$ lifts to a map which is the identity on the first component and which is constantly equal to the unit quaternion on the second one. In particular,
 \[\text{deg}(\eta_1)-\text{deg}(\eta_2)=1-0=1=\chi(\eta)\] and
 \[\text{deg}(\eta_1) + \text{deg}(\eta_2)=1+0=1=-p_1(\eta)/2\]
 which concludes the argument.
\end{proof}

\begin{remark}
    One can also prove equation (\ref{1}) by noticing that
    \begin{equation}
        \chi(M)=\text{deg}(f_1 \cdot f_2^{-1}).
    \end{equation}

    Indeed, recall that the double cover map
    \[\pi: S^3 \times S^3 \to SO(3)\]
    is defined by \[\pi(q_1,q_2)(q)\defeq q_1 \cdot q \cdot q_2^{-1}.\]
    This implies that the map 
    \[f_1 \cdot f_2^{-1}: S^3 \to S^3\]
    is recording the coordinates of a non-vanishing vector field on $M \setminus h^4$ with respect to the standard framing of the unit 4-ball $h^4$. The conclusion then follows from Poincaré-Hopf's theorem (see e.g. \cite[Chapter 6]{Milnor}) and from the fact that \[\text{deg}(f_1 \cdot f_2^{-1})=\text{deg}(f_1)-\text{deg}(f_2).\]
\end{remark}

\section*{Acknowledgements} The author is grateful to Younes Benyahia, Oliviero Malech, Rafael Torres and Daniele Zuddas for always supporting and encouraging her work and for their comments on this short note. She is also grateful to Marc Kegel for the nice discussions about this topic. The author has been partially supported by GNSAGA – Istituto Nazionale di Alta Matematica ‘Francesco Severi’, Italy.

\end{document}